# The Topological Directional Entropy of $\mathbb{Z}^2$-actions Generated by Linear Cellular Automata

## Hasan Akın

Department of Mathematics,
Arts and Sciences Faculty,
Harran University; 63100, Şanlıurfa, Turkey
E-mail address: akinhasan@harran.edu.tr

---

**Abstract:**

In this paper we study the topological directional entropy of $\mathbb{Z}^2$-actions by generated linear cellular automata (CA hereafter), defined by a local rule $f[l,r]$, $l, r \in \mathbb{Z}$, $l \leq r$, i.e. the maps $T_{f[l, r]}: \mathbb{Z}_m^{\mathbb{Z}} \to \mathbb{Z}_m^{\mathbb{Z}}$ which are given by $T_{f[l, r]}(x) = (y_n)_{n=-\infty}^{\infty}$, $y_n = f(x_{n+l}, ..., x_{n+r}) = \sum\limits_{i=l}^{r} \lambda_i x_{n+i} \pmod{m}$, $x = (x_n)_{n=-\infty}^{\infty} \in \mathbb{Z}_m^{\mathbb{Z}}$ and $f: \mathbb{Z}_m^{r-l+1} \to \mathbb{Z}_m$, over the ring $\mathbb{Z}_m$ ($m \geq 2$), and the shift map acting on compact metric space $\mathbb{Z}_m^{\mathbb{Z}}$. We give a closed formula, which can be efficiently and rightly computed by means of the coefficients of the local rule $f$, for the topological directional entropy of $\mathbb{Z}^2$-action generated by the pair $(T_{f[l, r]}, \sigma)$ in the direction $\theta$. We generalize the results obtained by Akın [The topological entropy of invertible cellular automata, J. Comput. Appl. Math. 213 (2) (2008) 501 − 508] to the topological entropy of any invertible linear CA.

Keywords: Invertible Linear Cellular Automata, Topological Entropy, Topological Directional Entropy.

---

## 1. Introduction

Cellular automata (CA), initialed by Ulam and von Neumann, has been systematically studied by Hedlund from purely mathematical point of view [12]. Hedlund's paper initiated investigation of current problems in symbolic dynamics. In [17], Shereshevsky has studied ergodic properties of CA, he has also defined $n$th iteration of permutative cellular automata and shown that if the local rule $f$ is right (left) permutative, then its $n$th iteration also is right (left) permutative. In [8], it has been answered open questions about the topological and ergodic dynamics of one-dimensional CA. Manzini and Margara [11] have obtained some necessary and sufficient conditions to be invertible for CA over $\mathbb{Z}_m$.

It is well known that there are several notions of entropy (i.e. measure entropy, topological entropy, directional entropy etc.) of measure-preserving transformation on probability space in ergodic theory. It is important to know how these notions are related with each other. In the last decade (see [1], [2], [4-10] and [14-15]), a lot papers are devoted to this subject.

Although the linear CA theory and the entropy of these linear CA have grown up somewhat independently, there are strong connections between entropy theory and CA theory. For the definitions and some properties of linear one-dimensional CA we refer to [3] and [14].

In [1], the author has proved that the uniform Bernoulli measure is a measure of maximal entropy for linear one-dimensional CA, he has studied the measure theoretical entropy of linear one-dimensional CA. In [11], D'amico et al. have studied the computing the topological entropy of $D$-dimensional CA by using an algorithm and Lyupanov exponents. They have given a closed formula for the Lyupanov exponents of 1-dimensional linear and positively expansive CA over the ring $\mathbb{Z}_m$ ($m \geq 2$). They have also shown in their paper that the topological entropy of CA is equal to the sum of the left and right Lyapunov exponents multiplied by the entropy of the shift CA provided that when the following property is satisfied: topological entropy obeys mentioned above if and only if the local associated with the CA is leftmost and rightmost permutative.

It is well known that if the finite *fps* associated with the local rule is invertible, so is the linear CA corresponding to this *fps*. Due to this fact, in order to study the invertibility of CA we consider the *fps*. The technique of *fps* is well known for the computation of the topological entropy and the finding of inverse of CA over the ring $\mathbb{Z}_m$ (see [4], [6] and [11] for details).

In general the calculation of the directional entropy of any $\mathbb{Z}^2$-action is difficult. Even, it is not almost possible. The concept of the directional entropy of a $\mathbb{Z}^2$-action has first been introduced by Milnor [14]. Milnor has defined the concept of directional entropy function for $\mathbb{Z}^2$-action generated by a full shift and a block map. Recently, the concept of the directional entropy which was initialed by Milnor has been studied by workers ([2], [8], [10] and [16]). In [15], Park has studied the continuity of directional entropy for arbitrary $\mathbb{Z}^2$-actions. In [10], it has been estimated from above the directional entropy of any CA-action with respect to an invariant Borel probability measure. They have considered permutative CA-actions.

In this paper, our main aim is to give an algorithm for computing the topological directional entropy of the $\mathbb{Z}^2$-actions generated by any 1-dimensional linear CA and the shift map. Thus, we ask to give a closed formula for the topological directional entropy of $\mathbb{Z}^2$-action generated by the pair $(T_{f[-r, r]}, \sigma)$ in the direction $\theta$ that can be efficiently and rightly computed by means of the coefficients of the local rule $f$ as similar to [11]. Also, we study the topological entropy of one-dimensional an invertible linear CA over $\mathbb{Z}_m$ with local rule $f$ by means of algorithm defined by D'amico et al. [11]. In this paper we are interested in 1-dimensional invertible CA and its topological entropy. We generalize the results obtained in [6] to any invertible CA over the ring $\mathbb{Z}_m$.

The rest of this paper is organized as follows: In Section 2, it is given basic definitions and notations. In Section 3, it is studied the topological entropies of any invertible one-dimensional linear CA and its inverses. In Section 4, it is given an algorithm for computing the topological directional entropy of the $\mathbb{Z}^2$-actions. In Section 4, we conclude by pointing some further problems.

## 2. Preliminaries

Let $\mathbb{Z}_m$ ($m \geq 2$) be a ring and $\mathbb{Z}_m^{\mathbb{Z}}$ be the space of all doubly-infinite sequences $x = (x_n)_{n=-\infty}^{\infty}$, $x_n \in \mathbb{Z}_m$. The shift $\sigma$: $\mathbb{Z}_m^{\mathbb{Z}} \to \mathbb{Z}_m^{\mathbb{Z}}$ defined by $(\sigma x)_i = x_{i+1}$ is a homeomorphism of compact metric space $\mathbb{Z}_m^{\mathbb{Z}}$. A CA is a continuous map, which commutes with $\sigma$, $T$: $\mathbb{Z}_m^{\mathbb{Z}} \to \mathbb{Z}_m^{\mathbb{Z}}$ defined for $x = (x_n)_{n=-\infty}^{\infty} \in \mathbb{Z}_m^{\mathbb{Z}}$, $i \in \mathbb{Z}$ by $(Tx)_i = f(x_{i-r}, ...., x_{i+r})$, where $f$: $\mathbb{Z}_m^{2r+1} \to \mathbb{Z}_m$ is a given local rule or map. The integer $r$ is called the radius of the CA. It has been defined that the local rule $f$ is linear if and only if it can be written as $f(x_l, ..., x_r) = \sum_{i=l}^{r} \lambda_i x_i \pmod{m}$, where at least one between $\lambda_l$ and $\lambda_r$ is nonzero. We consider linear one-dimensional CA $T_{f[l, r]}$ determined by local rule $f$:

$$Tx = (\ y_n\ )_{n=-\infty}^{\infty}\ ,\ \ y_n = f(x_{n+l}, ..., x_{n+r}) = \sum_{i=l}^{r} \lambda_i x_{n+i} \pmod{m}, \qquad (2.1)$$

where $\lambda_l, .., \lambda_r \in \mathbb{Z}_m$. The notation of permutive CA was first introduced by Hedlund in [12]. Afterwards, in [17] Shereshevsky has defined the permutative notation in a different way as follows; given a mapping $f:\ \mathbb{Z}_m^r \rightarrow \mathbb{Z}_m$ and a block $(\hat{x}_1, ..., \hat{x}_{r-1}) \in \mathbb{Z}_m^{r-1}$ defined a map

$$f_{(\hat{x}_1, ..., \hat{x}_{r-1})}^{+}(x_r) = f(\hat{x}_1, ..., \hat{x}_{r-1}, x_r) = (\sum_{i=1}^{r-1} \lambda_i \hat{x}_i + \lambda_r x_r)(\mathrm{mod}\, m)\,.$$

Likewise, for any $(\hat{x}_2, ..., \hat{x}_r) \in \mathbb{Z}_m^{r-1}$ it is shown

$$f_{(\hat{x}_2, ..., \hat{x}_r)}^{+}(\ x_1) = f(x_1, \hat{x}_2, ..., \hat{x}_r) = (\sum_{i=2}^{r} \lambda_i \hat{x}_i + \lambda_1 x_1)(\mathrm{mod}\, m)\,.$$

It is said that the rule is right (resp., left) permutative, if the mapping $f_{(\hat{x}_1, ..., \hat{x}_{r-1})}^{+}$ (resp., $f_{(\hat{x}_2, ..., \hat{x}_r)}^{+}$) is a permutation of $\mathbb{Z}_m$ for every $(\hat{x}_1, ..., \hat{x}_{r-1}) \in \mathbb{Z}_m^{r-1}$ (resp., $(\hat{x}_2, ..., \hat{x}_r) \in \mathbb{Z}_m^{r-1}$).

Note that $f$ is permutative in the $j$th variable if and only if gcd $(\lambda_j, m) = 1$. In this paper we will use this definition.

The mapping $f$ is called the (generating) local rule of the CA-map. We will use the notation $T_{f[l, r]}$ for linear CA-map defined (2.1) to emphasize the local rule $f$ and the numbers $l$ and $r$. If the local rule $f$ is as $f(x_l, ..., x_r) = \sum_{i=l}^{r} \lambda_i x_i \pmod{m}$, then the finite formal power series (*fps* for short) $F$ associated with $f$ is $F(X) = \sum_{i=l}^{r} \lambda_i X^i$. We will use notations the local rule $f$, generic finite formal power series (*fps*) $F$ and 1-dimensional CA $T_{f[l, r]}$ ( see [11] and [13] for details).

### 3. Topological Entropy of Invertible CA

The entropy has been interpreted as a measure of the chaotic character of a dynamical system by very much authors (see [8]), the value $\mathcal{H}(\mathbb{Z}_m^{\mathbb{Z}}, T)$ has been in general accepted as a measure of the complexity of the dynamics of $T$ over $\mathbb{Z}_m^{\mathbb{Z}}$ [16] or the topological entropy of a map is a crude global measure of the exponential complexity of the structure of the orbits of the map. Also it has been stated that the topological entropy of an ordinary dynamical system is, roughly speaking, the measure of the exponential increase with time of the maximal number of distinct trajectories up to some precision.

In [1] Akın has calculated the measure entropy of linear one-dimensional CA with respect to uniform Bernoulli measure. Recall that by the Variational Principle the topological entropy is the supremum of the entropies of invariant measures. In [1], it has been shown that the uniform Bernoulli measure is a measure of maximal entropy for linear one-dimensional CA. Thus, these measures carry the maximal amount of information of the system. We conclude that the value $\mathcal{H}(\mathbb{Z}_m^{\mathbb{Z}}, T_{f[l, r]})$ is the maximal amount of information that can be carried by the uniform Bernoulli measure (see [6]). See [11] for the definition of computing the topological entropy calculated by means of the algorithm and the Lyapunov exponents. Now we can state the algorithm given by [11].

**Theorem 3.1** ([11], Theorem 2). Let $(\mathbb{Z}_m^{\mathbb{Z}}, T_{f[-r, r]})$ be a 1-dimensional CA over $\mathbb{Z}_m$ with local rule $f(x_{-r}, \ldots, x_r) = \sum\limits_{i=-r}^{r} \lambda_i x_i \pmod{m}$ and let $m = p_1^{k_1} \ldots p_h^{k_h}$ be the prime factor decomposition of m. For $i = 1,\ldots, h$ define $P_i = \{0\} \cup \{j : \gcd(\lambda_j, p_i) = 1\}$, $L_i = \min P_i$, $R_i = \max P_i$. Under the properties of theorem the topological entropy of $(\mathbb{Z}_m^{\mathbb{Z}}, T_{f[-r, r]})$ is

$$\mathcal{H}(\mathbb{Z}_m^{\mathbb{Z}}, T_{f[-r, r]}) = \sum_{i=1}^{h} k_i (R_i - L_i) \log(p_i).$$

In [13], Manzini and Margara have studied the problem of finding the inverse of a linear $D$-dimensional CA over $\mathbb{Z}_m$. To find the inverse of a linear CA is difficult because the radius of $T^{-1}$ is in general different from the radius of $T$. This fact is given in [13]. A CA is invertible if its global transformation map which is obtained by applying the local rule $f$ to all sites of the lattice is invertible.

**Definition 3.2.** If $T_{f[l, r]}: \mathbb{Z}_m^{\mathbb{Z}} \to \mathbb{Z}_m^{\mathbb{Z}}$ is the linear CA determined by local rule $f$:

$$(T_{f[l, r]} x)_n = \sum_{i=l}^{r} \lambda_i x_{n+i} \pmod{m},$$

then $T_{f[l, r]}$ is invertible if and only if for each prime factor $p$ of $m$ there exists a unique coefficient $\lambda_j$ ($l \leq j \leq r$) such that $\gcd(p, \lambda_j) = 1$, that is, $p \mid \lambda_j$ and $\gcd(p, \lambda_j) = 1$ for $i \neq j$.

We can state the following Lemma for 1-dimensional CA.

**Lemma 3.3** ([13], Lemma 3.1). Let $m = p_1^{k_1} \cdot p_2^{k_2} \ldots p_h^{k_h}$ and $F(X)$ be a *fps* over $\mathbb{Z}_m$. Given $h$ finite fps's $G_1, \ldots, G_h$ such that $F(X).G_i(X) \equiv 1 \pmod{P_i^{k_i}}$, it can be found a h finite fps $G$ such that $F(X).G(X) \equiv 1 \pmod{m}$. So it is found $G(X) = \sum\limits_{i=1}^{h} \beta_i \alpha_i G_i(X)$, where $\alpha_i = \dfrac{m}{p_i^{k_i}}$ and $\alpha_i \beta_i \equiv 1 \pmod{p_i^{k_i}}$. $F$ can be written as

$$F(X) = \lambda_j X^j + p H(X),$$

where $\gcd(\lambda_i, p) = 1$.

**Theorem 3.4** ([13], Theorem 3.2). Let $F(X)$ denote an invertible finite *fps* over $\mathbb{Z}_{p^k}$, and let $\lambda_j$ and $H$ be defined as in (5). Let $\lambda_j^{-1}$ be such that $\lambda_j^{-1} . \lambda_j \equiv 1 \pmod{p_i^{k_i}}$. Then, the inverse of $F$ is given by $G(X) = \lambda_j^{-1} X^{-j} (1 + p \tilde{H}(X) + p^2 \tilde{H}^2(X) + \ldots + p^{k-1} \tilde{H}^{k-1}(X))$, where $\tilde{H}(X) = -\lambda_j^{-1} X^{-j} H(X)$.

In this section, our main result is the following Theorem.

**Theorem 3.5.** Let $m = p_1^{k_1} \cdot p_2^{k_2} \dots p_h^{k_h}$ and $F(X) = \sum_{i=1}^{r} \alpha_i X^{-i}$ be an *fps* over $\mathbb{Z}_m$. Given $h$ finite *fps*'s $G_1, \dots, G_h$ such that $F(X).G_i(X) \equiv 1 \pmod{P_i^{k_i}}$. If the local rule $f$ is associated with the *fps* $F$, then the topological entropy of CA generated by the local rule $f$ is as follows; $\mathcal{H}\ell(\mathbb{Z}_m^{\mathbb{Z}}, T_{f[l,\, r]}) = \sum_{i=1}^{h} j_i \log p_i^{k_i}$, where for all $i = 1, 2, \dots, h$ $\alpha_{j_i} = \dfrac{m}{p_i^{k_i}}$, $gcd(\alpha_{j_i}, p_i) = 1$.

From Theorem 3.1, the proof follows.

It is easy to see that $\mathcal{H}(\mathbb{Z}_m^{\mathbb{Z}}, T_{f[l,\, r]}) = \mathcal{H}(\mathbb{Z}_m^{\mathbb{Z}}, T_{f[l,\, r]}^{-1})$, where $T_{f[l,\, r]}^{-1}$ is the inverse of invertible CA $T_{f[l,\, r]}$.

# 4. Directional Entropy.

In this section, our aim is to obtain an algorithm in order to compute efficiently the topological directional entropy as similar to the algorithm defined by D'amico at al. [11]. In [9], the expanded proof of the topological directional entropy has been given by the following formula (see [9] for details)

$$h(\theta) = \begin{cases} (\cos\theta + r\sin\theta)\ln p; & \theta \in [0, \theta_l] \\ (r-l)\sin\theta\ln p; & \theta \in [\theta_l, \theta_r] \\ (\cos\theta + l\sin\theta)\ln p; & \theta \in [\theta_r, \pi] \end{cases}, \qquad (4.1)$$

where $\theta_l = \operatorname{arccot}(-l)$ and $\theta_r = \operatorname{arccot}(-r)$, also the local rule $f$ defined in (2.1) is bipermutative.

Let us give the following Lemma. This Lemma is important to compute in general case the topological directional entropy of any $\mathbb{Z}^2$-action by the pair $(T_{f[-r,\, r]}, \sigma)$, where $T_{f[-r,\, r]}$ is any linear CA.

**Lemma 4.1.** Let $(\mathbb{Z}_{p^k}^{\mathbb{Z}}, T_{F[-r,\, r]})$ be the one-dimensional linear CA over $\mathbb{Z}_{p^k}$ with local rule $f(x_{-r}, \dots, x_r) = \sum_{i=-r}^{r} \lambda_i x_i \pmod{p^k}$, where $p$ is prime number. If

$$P = \{0\} \cup \{j: \gcd(\lambda_j, p)\} = \{j_1, j_2, \dots, j_t\}, L = \min P \text{ and } R = \max P. \qquad (4.2)$$

Then the topological directional entropy of $\mathbb{Z}^2$-action generated by the pair $(T_{f[-r,\, r]}, \sigma)$ in the direction $\theta$ is

$$h(\theta) = \begin{cases} k(\cos\theta + R\sin\theta)\ln p; & \theta \in [0, \theta_L] \\ k(R-L)\sin\theta\ln p; & \theta \in [\theta_L, \theta_R] \\ k(\cos\theta + L\sin\theta)\ln p; & \theta \in [\theta_R, \pi]. \end{cases} \qquad (4.3)$$

**Proof:** It is well known that $f$ has radius $r$, i.e., it depends on at most $2r+1$ variables. We associate to the local rule $f$ the finite *fps* $F(X) = \sum_{i=-r}^{r} \lambda_i X^{-i}$. It is well known that the finite *fps* associated with $f^{(n)}$ is $F^n(X)$. From (Lemma 1, [11]), it is easy to see that there exist $L, R \in \mathbb{Z}$ and $n \in \mathbb{N}$ such that $f^{(n)}$ is permutative in variables $x_L$ and $x_R$. Furthermore, the local rule $f^{(n)}$ does not depend on variables $x_j$ with $j<L$ or $j>R$; that is we can ignore the variables $x_j$ with $j<L$ or $j>R$. Then, the local rule $f^{(n)}$ is both *leftmost* and *rightmost* permutative. Thus, from (4.1) and (Lemma 1, [11]) we have (4.3). Thus, the proof is completed.

**Theorem 4.2.** Let $(\mathbb{Z}_m^{\mathbb{Z}}, T_{f[-r, r]})$ be a 1-dimensional CA over $\mathbb{Z}_m$ with local rule

$$f(x_{-r}, ..., x_r) = \sum_{i=-r}^{r} \lambda_i x_{n+i} \pmod{m}$$ and let $m = p_1^{k_1} ... p_h^{k_h}$ be the prime factor decomposition of $m$. For $i = 1, ..., h$ define $P_i = \{0\} \cup \{j: \gcd(\lambda_j, p_i)=1\}$, $L_i = \min P_i$, $R_i = \max P_i$. Then the topological directional entropy of $\mathbb{Z}^2$-action generated by the pair $(T_{f[-r, r]}, \sigma)$ in the direction $\theta$ is the following

$$h(\theta) = \begin{cases} \sum_{i=1}^{h} k_i(\cos\theta + R_i\sin\theta)\ln p_i; & \theta \in [0, \min\theta_{L_j}]; \\[2mm] \sum_{\substack{i=1 \\ i \neq j}}^{h} k_i(\cos\theta + R_i\sin\theta)\ln p_i & \theta \in [\min\theta_{L_j}, \min\{\theta_{Rj_1}, \theta_{L_{j_1}}\}]; \\[1mm] + k_j(R_j - L_j)\sin\theta\ln p; & \\[1mm] \quad\ldots\ldots\ldots\ldots & \quad\ldots\ldots \\[1mm] \sum_{i=1}^{h} k_i(\cos\theta + L_i\sin\theta)\ln p_i; & \theta \in [\max\theta_{R_j}, \pi]. \end{cases} \quad (4.4)$$

Proof: The proof easily can be obtained from Lemma 4.1.

**Example 4.3.** Let us consider the local rule
$$f(x_{-3}, ..., x_2) = (15x_{-3} + 20x_{-2} + 27x_{-1} + 16x_0 + 30x_1 + 5x_2) \pmod{2^3 3^5 5^2}$$

and let $T_{f[-3, 2]}$ be the 1-dimensional linear CA associated with $f$. From Theorem 3.1, that is, from the formula we can obtain as $P_{p_i} = \{0\} \cup \{j \mid \gcd(\lambda_j, p_i) = 1\}$. Thus, we can write as $P_2 = \{-3, 0, 2\}$, $P_3 = \{-2, 0, 2\}$, $P_5 = \{-1, 0\}$. For primes $p_i = 2, 3, 5$ define $\mathbb{P}_2 = \{\cot\theta_2 = -2, \cot\theta_{-3} = 3\}$, $\mathbb{P}_3 = \{\cot\theta_2 = -2, \cot\theta_{-2} = 2\}$, $\mathbb{P}_5 = \{\cot\theta_0 = 0, \cot\theta_{-1} = 1\}$. Then from 4.2, we have

$$h(\theta) = \begin{cases} \cos\theta\ln 2^3 3^5 5^2 + \sin\theta\ln 2^6 3^{10}, & \theta \in [0, \theta_{-3}]; \\[2mm] \cos\theta\ln 3^5 5^2 + \sin\theta\ln 2^{15} 3^{10} 5^4, & \theta \in [\theta_{-3}, \theta_{-2}]; \\[2mm] \cos\theta\ln 5^2 + \sin\theta\ln 2^{15} 3^{20}, & \theta \in [\theta_{-2}, \frac{\pi}{4}]; \\[2mm] \sin\theta\ln 2^{15} 3^{20} 5^2, & \theta \in [\frac{\pi}{4}, \frac{\pi}{2}]; \\[2mm] \cos\theta\ln 5^2 + \sin\theta\ln\dfrac{2^{15} 3^{20}}{5^2}, & \theta \in [\frac{\pi}{2}, \theta_2]; \\[2mm] \cos\theta\ln 2^3 3^5 5^2 - \sin\theta\ln 2^9 3^{10} 5^2, & \theta \in [\theta_2, \pi]. \end{cases}$$

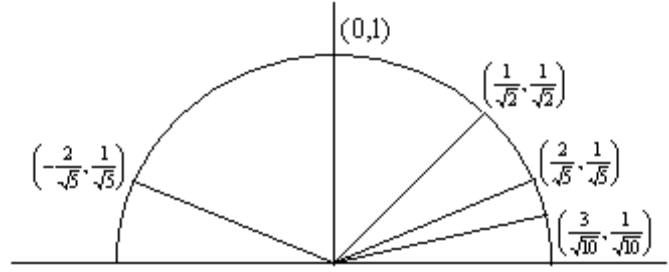

**Figure 1**

From Figure 1, we can easily see that for every sector the value $h(\theta)$ is different.

Now let us study the topological directional entropy of $\mathbb{Z}^2$-action generated by the pair $(T_{f[l, r]}, \sigma)$, where $T_{f[l, r]}$ is invertible CA.

**Proposition 4.4.** Let local rule $f$ be as follows;

$$f(x_l, \ldots, x_r) = \lambda_j x_j + p \sum_{\substack{i=l \\ i \neq j}}^{r} a_i x_i \ (\mathrm{mod} \ p^k),\tag{4.5}$$

where $p \nmid \lambda_j$. Then the topological directional entropy of $\mathbb{Z}^2$-action generated by the pair $(T_{f[-r, r]}, \sigma)$ in the direction $\theta$ for any $\theta \in [0, \pi]$

$$h(\theta) = k(\cos\theta + j\sin\theta)\ln p.\tag{4.6}$$

**Proof:** It is clear that for the local rule $f$ we have $P = \{0\} \cup \{j: \gcd(\lambda_j, p) = 1\} = \{0, j\}$. Thus from Lemma 4.1 we have $h(\theta) = k(\cos\theta + j\sin\theta)\ln p$. We can conclude that if $\theta = \dfrac{\pi}{2}$, then the topological directional entropy reaches the maximal value.

Let us consider the local rule $f$ in (4.5). If $F(X)$ is the finite *fps* associated with the CA generated by the local rule $f$, then we can write as follows;

$$F(\mathrm{X}) = \lambda_j X^{-j} + p \sum_{\substack{i=l \\ i \neq j}}^{r} a_i X^{-i} = \lambda_j X^{-j} + H(X).\tag{4.7}$$

Then the inverse of $F$ is given by

$$G(\mathrm{X}) = \lambda_j^{-1} \mathrm{X}^{j}(1 + p\widetilde{H} + p^2 \widetilde{H}^2 + \ldots + p^{k-1} \widetilde{H}^{k-1}),\tag{4.8}$$

where $\widetilde{H}(\mathrm{X}) = -\lambda_j^{-1} \mathrm{X}^{j} H(X)$ and $p \nmid \lambda_j$. From (4.8) we can write the local rule associated with the *fps* G(X) as follows;

$$g(\mathrm{x}_{l(k-1)+j}, \ldots, \mathrm{x}_{r(k-1)+j}) = \lambda_j \ x_{-j} + p \sum_{\substack{i=l(k-1)+j \\ i \neq -j}}^{l(k-1)+j} a_i x_i \ (\mathrm{mod} \ p^k).\tag{4.9}$$

So from Lemma 4.1 we conclude that the topological directional entropy of $\mathbb{Z}^2$-action generated by $T_{g[l(k-1)+j, \ r(k-1)+j]}$ and the shift map in direction $\theta$ is

$$h(\theta) = k(\cos\theta + j\sin\theta)\ln p. \tag{4.10}$$

Thus from (4.6) and (4.10) we conclude that if $\theta = \dfrac{\pi}{2}$, then the directional entropy coincides with the topological entropy of the $\mathbb{Z}^2$-action generated by both $(T_{g[l(k-1)+j,\ r(k-1)+j]}\ ,\ \sigma)$ and $(T_{f[l,\ r]}\ ,\ \sigma)$.

**Example 4.5.** Let us consider the local rule $f(x_1, x_2, x_3, x_4) = 2(x_1 + x_2 + x_3) + x_4 (\mathrm{mod}\ 2^2)$.

Then we have $\cot\theta_l = -4$ and $\cot\theta_r = -4$, thus the topological directional entropy of the $\mathbb{Z}^2$-actions generated by the pair $(T_{f[1,\ 4]}\ ,\ \sigma)$ is $h(\theta) = 2(\cos\theta + 4\sin\theta)\ln 2$, where $\theta \in [0, \pi]$. The inverse of $f$ is $g(x_{-4}, x_{-6}, x_{-5}, x_{-4}) = 2(x_{-7} + x_{-6} + x_{-5}) + x_{-4} (\mathrm{mod}\ 2^2)$, so we have $\cot\theta_l = 4$, $\cot\theta_r = 4$ for all $\theta \in [0, \pi]$. Thus from Proposition 4.3, the topological directional entropy of the $\mathbb{Z}^2$-actions generated by the pair $(T_{g[-7,\ -4]}\ ,\ \sigma)$ is $h(\theta) = 2(\cos\theta + 4\sin\theta)\ln 2$.

# 5. Conclusions

In this paper, we showed how to compute the topological entropy of any invertible CA and the topological directional entropy of $\mathbb{Z}^2$-action generated by CA as similar to the algorithm defined by D'amico at al. [11] (see also [9]). These results are generalizations of previous work done by several authors [6] and [9-11]. It is also possible to extend these results to a finite product for Galois rings (see [5]). Similar computations and explorations of CA over different rings remain to be of interest.